\documentclass[12pt]{amsart}
\usepackage{amscd,amsmath,amssymb,amsfonts}
\usepackage[all]{xy}
\theoremstyle{plain}
\newtheorem{thm}{Theorem}
\newtheorem{lem}[thm]{Lemma}

\theoremstyle{definition}
\newtheorem{defn}[thm]{Definition}

\numberwithin{thm}{section} \numberwithin{equation}{section}

\newcommand{\sO}{{\mathcal O}}

\begin{document}

\title{stable bundles as Frobenius morphism direct image}
\author{Congjun Liu}
\address{Academy of Mathematics and Systems Science, Chinese Academy of Science, Beijing, P. R. of China}
\email{liucongjun@amss.ac.cn}

\author{Mingshuo Zhou}
\address{Academy of Mathematics and Systems Science, Chinese Academy of Science, Beijing, P. R. of China}
\email{zhoumingshuo@amss.ac.cn}
\date{November 29, 2012}

\begin{abstract}
Let $X$ be a smooth projective curve of genus $g\geq 2$  over an
algebraically closed field $k$ of characteristic $p>0$, and let
$F:X\rightarrow X_{1}$ be the relative  Frobenius map. We show that
a vector bundle E on $X_{1}$ is the direct image of some stable
bundle $X$ if and only if instability of $F^{*}E$ is equal to $(p-1)(2g-2)$.
\end{abstract}
\maketitle

\section{Introduction}
Let $X$ be a smooth projective curve of genus $g\geq 2$ defined over
an algebraically closed field $k$ of characteristic $p>0$. The absolute
Frobenius morphism $F_{X}$ : $X\rightarrow X$ is induced by
$\sO_{X} \rightarrow \sO_{X}, f\mapsto f^p$. Let $F: X \rightarrow X_{1}:=X\times_{k}k$
denote the relative Frobenius morphism over $k$. One of the themes is to study
its action on the geometric objects on $X$.
Recall that a vector bundle $E$ on a smooth projective curve is called semi-stable (resp. stable)
if $\mu(E') \leq \mu(E)$ (resp. $\mu(E') < \mu(E)$) for any nontrivial proper subbundle $E' \subset E$,
where $\mu(E)$ is the slope of $E$. It is known that $F_{*}$ is preserves
the stability of vector bundles (cf.\cite{Sun1}), but $F^{*}$ does not preserve
the semi-stability of vector bundle(cf.\cite{G} for example).

Semi-stable bundles are basic constituents
of vector bundles in the sense that any bundles $E$ admits a unique filtration
$${\rm HN}_{\bullet}(E): 0={\rm HN}_0(E)\subset{\rm HN}_1(E)\subset \cdots\subset{\rm HN}_{\ell}(E)=E,$$
  which is the so called Harder-Narasimhan filtration, such that
  \begin{itemize}
 \item[\rm{(1)}] ${\rm gr}_i^{\rm HN}(E):={\rm HN}_i(E)/{\rm HN}_{i-1}(E)$ ($1\le i\le \ell$) are semistable;
\item[\rm{(2)}]  $\mu({\rm gr}_1^{\rm HN}(E))>\mu({\rm gr}_2^{\rm HN}(E))>\cdots>\mu({\rm gr}_{\ell}^{\rm HN}(E))$.
  \end{itemize}
The rational number ${\rm I}(E) :=\mu({\rm gr}_1^{\rm
HN}(E))-\mu({\rm gr}_{\ell}^{\rm HN}(E))$, which measures how
far is a vector bundle from being semi-stable, is called the
instability of $E$. It is clear that $E$ is semi-stable if and
only if ${\rm I}(E)=0$.

Given a semi-stable bundle $E$ on $X_{1}$, then $F^{*}E$ may not be semi-stable, so it is natural to consider the instability
${\rm I}(F^{*}E)$. In \cite[Theorem~3.1]{Sun2}, the author prove ${\rm I}(F^{*}E) \leq (\ell-1)(2g-2)$,
where $\ell$ is the length of  Harder-Narasimhan filtration of $F^{*}E$. If $E=F_*W$ where $W$
is stable bundle on $X$, we know, by Sun's theorem (\cite [theorem~2.2]{Sun1}), that $E$ is stable, the length
of Harder-Narasimhan filtration of $F^*E$ is $p$ and ${\rm I}(F^{*}E)=(p-1)(2g-2)$. Thus ${\rm I}(F^{*}E)=(p-1)(2g-2)$
is a necessary condition that $E$ is a direct image under Frobenius.
In this short note, we show the following theorem:
\begin{thm}Let E be a stable vector bundle over X. Then
the following statements are equivalent:
\item{(1)} There exists a stable bundle W such that $E=F_{*}W$;
\item{(2)} ${\rm I}(F^{*}E)=(p-1)(2g-2).$
\end{thm}
\noindent The case ${\rm rk}E=p$ was proved in \cite{MP}. Our observation is that the arguments
in \cite{MP} and Sun's theorem together imply the general case.

\section{Proof of the theorem}
Let $X$ be a smooth projective curve over an
algebraically closed field $k$ with ${\rm char}(k)=p>0$. The
absolute Frobenius morphism $F_X:X\to X$ is induced by the
homomorphism $$\mathcal{O}_X\to \mathcal{O}_X,\qquad f\mapsto f^p$$
of rings. Let $F:X\to X_1:=X\times_kk$ denote the relative Frobenius
morphism over $k$ that satisfies
$$\xymatrix{\ar@/^20pt/[rr]^{F_X} X\ar[r]^F
\ar[dr] & X_1\ar[r]\ar[d]
& X\ar[d]^{}\\
 & {\rm Spec}(k)\ar[r]^{F_k}& {\rm Spec}(k)} .$$

For a vector bundle $E$ on $X$,the slope of $E$ is defined as
$$\mu(E):= \frac{{\rm deg} E}{{\rm rk} E}$$ where ${\rm rk} E$ (resp. deg $E$) denotes the rank (resp. degree) of $E$.
 Then
\begin{defn} A vector bundle $E$ on a $X$ is called semi-stable (resp. stable) if for any nontrivial proper
subbundle $E' \subset E$, we have $$\mu(E') \leq (resp. <) \mu(E).$$
\end{defn}

\begin{thm}(Harder-Narasimhan filtration)
  For any vector bundle $E$, there is a unique filtration
  $${\rm HN}_{\bullet}(E): 0={\rm HN}_0(E)\subset{\rm HN}_1(E)\subset \cdots\subset{\rm HN}_{\ell}(E)=E,$$
  which is the so called Harder-Narasimhan filtration, such that
  \begin{itemize}
 \item[\rm{(1)}] ${\rm gr}_i^{\rm HN}(E):={\rm HN}_i(E)/{\rm HN}_{i-1}(E)$ ($1\le i\le \ell$) are semistable;
\item[\rm{(2)}]  $\mu({\rm gr}_1^{\rm HN}(E))>\mu({\rm gr}_2^{\rm HN}(E))>\cdots>\mu({\rm gr}_{\ell}^{\rm HN}(E))$.
  \end{itemize}
  \end{thm}

By using this unique filtration of $E$, an
invariant ${\rm I}(E)$ of $E$, which is called the instability of
$E$ was introduced (see \cite{Sun1} and \cite{Sun2}). It is a rational number and measures how far is $E$ from
being semi-stable.

\begin{defn} Let $\mu_{\rm max}(E)=\mu({\rm gr}_1^{\rm
HN}(E))$,
  $\mu_{\rm min}(E)=\mu({\rm gr}_{\ell}^{\rm HN}(E))$. Then the instability of $E$ is
  defined to be
  $${\rm I}(E):=\mu_{\rm max}(E)-\mu_{\rm min}(E).$$
\end{defn}

\noindent It is easy to see that a torsion free sheaf $E$ is
semi-stable if and only if ${\rm I}(E)=0$.

For any semi-stable bundle $E$, let
$${\rm HN}_{\bullet}(F^{*}E): 0={\rm HN}_0(F^{*}E)\subset{\rm HN}_1(F^{*}E)\subset \cdots\subset{\rm HN}_{\ell}(F^{*}E)=F^*E$$
be the Harder-Narasimhan filtration of $F^{*}E$. Then we have the following lemma, which is implicit in \cite{MP}.
\begin{lem} For any semi-stable bundle $E$, we have
$$\mu_{max}(F^*E)\leq p\cdot \mu(E)+(p-1)(g-1);$$$$\mu_{min}(F^*E)\geq p\cdot \mu(E)-(p-1)(g-1),$$
and if ${\rm I}(F^*E)=\mu_{\rm max}(F^*E)-\mu_{\rm min}(F^*E)=(p-1)(2g-2).$ Then
$$\mu_{max}(F^*E)= p\cdot \mu(E)+(p-1)(g-1);$$$$\mu_{min}(F^*E)= p\cdot \mu(E)-(p-1)(g-1).$$
\end{lem}

Now we prove our theorem by using this lemma and Sun's theorem on stability of Frobenius direct images.
\begin{proof} [Proof of Theorem 1.1]
$(1)\Rightarrow (2)$ is contained in \cite{Sun1}.\\
We prove $(2)\Rightarrow (1)$ here. Since ${\rm I}(F^*E)=(p-1)(2g-2)$, we have $\mu_{max}(F^*E)= p\cdot \mu(E)+(p-1)(g-1),\mu_{min}(F^*E)= p\cdot \mu(E)-(p-1)(g-1)$
by lemma 2.4. We consider the surjection $$F^*E\rightarrow gr_{\ell}^{HN}(F^*E).$$
The bundle $gr_{\ell}^{HN}(F^*E)$ is semi-stable of slope $\mu_{min}(F^*E)$. Replaced $gr_{\ell}^{HN}(F^*E)$ by a stable
graded piece $W$ in Jordan-H$\ddot{o}$lder filtration of $gr_{\ell}^{HN}(F^*E)$, we have a surjection
$$F^*E\rightarrow W,$$ where $W$ is a stable bundle of slope $\mu(W) = \mu_{min}(F^*E)= p\cdot \mu(E)-(p-1)(g-1).$
By adjunction we have a non-trivial morphism $$\psi: E\rightarrow F_*W.$$
By Sun's theorem (cf. \cite[Theorem~2.2]{Sun1}), we know that $F_*W$ is a stable bundle of slope
$$\mu(F_*W)= \frac{\mu(W)}{p}+\frac{(p-1)(g-1)}{p}= \mu (E).$$
Thus $\psi$ induce an isomorphism:
$$E\cong F_{*}W.$$
\end{proof}

\noindent 
\begin{itshape}Acknowledgements: The authors would like to thank their advisor professor Xiaotao Sun for encouragements 
and many useful discussions.
\end{itshape}

\bibliographystyle{plain}

\renewcommand\refname{References}

\end{document}